\documentclass[11pt]{article}
\usepackage[a4paper]{anysize}\marginsize{3.5cm}{3.5cm}{1.3cm}{2cm}
\pdfpagewidth=\paperwidth \pdfpageheight=\paperheight
\usepackage{amsfonts,amssymb,amsthm,amsmath,eucal}
\usepackage{pgf}
\usepackage{bbm,array}
\theoremstyle{plain}
\newtheorem{thm}{Theorem}[section]
\newtheorem{theorem}[thm]{Theorem}

\newtheorem{proposition}[thm]{Proposition}

\theoremstyle{definition}
\newtheorem{definition}[thm]{Definition}
\newtheorem{remark}[thm]{Remark}
\newtheorem{example}[thm]{Example}

\newtheorem{thevarthm}[thm]{\varthmname}

\newenvironment{varthm*}[1]{\trivlist\item[]{\bf #1.}\it}{\endtrivlist}

\renewcommand\geq{\geqslant}

\renewcommand\leq{\leqslant}

\newcommand\be{\begin{eqnarray*}}
\newcommand\ee{\end{eqnarray*}}

\newcommand\R{\mathbb R}

\newcommand\K{\mathbb K}
\newcommand\T{\mathbb T}
\renewcommand\P{\mathbb P}

\newcommand\calo{{\mathcal O}}
\newcommand\cali{{\mathcal I}}

\newcommand\newop[2]{\def#1{\mathop{\rm #2}\nolimits}}
\newop\log{log}
\newop\ord{ord}
\newop\Gal{Gal}
\newop\SL{SL}
\newop\GL{GL}
\newop\Bl{Bl}
\newop\mult{mult}
\newop\mass{mass}
\newop\div{div}
\newop\codim{codim}
\newop\sing{sing}
\newop\vdim{vdim}
\newop\edim{edim}
\newop\Ass{Ass}
\newop\size{size}
\newop\reg{reg}
\newop\areg{areg}
\newop\asreg{asreg}
\newop\satdeg{satdeg}
\newop\supp{supp}
\newop\gin{gin}
\newop\ini{in}
\newop\vol{vol}
\newop\length{length}
\newcommand\eqnref[1]{(\ref{#1})}

\newcommand\ibul{I_{\bullet}}
\newcommand\gbul{G_{\bullet}}

\def\keywordname{{\bfseries Keywords}}%
\def\keywords#1{\par\addvspace\medskipamount{\rightskip=0pt plus1cm
\def\and{\ifhmode\unskip\nobreak\fi\ $\cdot$
}\noindent\keywordname\enspace\ignorespaces#1\par}}
\def\subclassname{{\bfseries Mathematics Subject Classification
(2000)}\enspace}
\def\subclass#1{\par\addvspace\medskipamount{\rightskip=0pt plus1cm
\def\and{\ifhmode\unskip\nobreak\fi\ $\cdot$
}\noindent\subclassname\ignorespaces#1\par}}

\begin{document}

\author{M.~Dumnicki, T.~Szemberg\footnote{The second named author was partially supported
by NCN grant UMO-2011/01/B/ST1/04875}, J.~Szpond, H.~Tutaj-Gasi\'nska}
\title{Symbolic generic initial systems of star configurations}
\date{\today}
\maketitle
\thispagestyle{empty}

\begin{abstract}
   The purpose of this note is to describe limiting shapes
   of symbolic generic initial systems of star configurations
   in projective spaces $\P^n$ over a field $\K$ of characteristic $0$.
\subclass{MSC 14C20 \and MSC 13C05 \and MSC 14N05 \and MSC 14H20 \and MSC 14A05}
\end{abstract}

%*****************************************************************************

\section{Introduction}
\label{intro}
   In recent years there has been increasing interest in asymptotic invariants
   attached to graded families of ideals, see e.g. \cite[Section 2.4.B]{PAG}.
   For a homogeneous ideal $I$, Mayes introduces in \cite{May14}
   \emph{symbolic generic initial systems} $\left\{\gin(I^{(m)})\right\}_m$.
   Here $\gin(J)$ denotes the reverse lexicographic generic initial ideal
   of a homogeneous ideal $J$ and $J^{(m)}$ denotes the $m$th symbolic
   power of $J$. For a monomial ideal $J$ one defines its Newton polytope $P(J)$.
   Mayes studies the \emph{limiting shape} associated to $I$ as the
   complement in the positive octant in $\R^n$ of the asymptotic Newton polytope
   $$\Delta(I)=\bigcup_{m=1}^{\infty}\frac{1}{m}P(\gin(I^{(m)}))$$
   for \emph{generic} points in $\P^2$ under assuming the
   Segre-Harbourne-Gimigliano-Hirschowitz Conjecture, (see \cite{RDLS}
   for recent account on this conjecture).
   In this note we study limiting shapes for star configurations
   in projective spaces of arbitrary dimension.
   Star configurations have received much attention recently,
   partly because they are a nice source of interesting
   and computable examples, see e.g. the nice survey \cite{GHM}.
   Our main result here is the following theorem.
\begin{theorem}\label{thm:main points}
   Let $I$ be the ideal of points defined as $n$--fold intersection
   points of $s\geq n$ general hyperplanes in $\P^n$. Then
   $\Gamma(I)=\overline{(\R_{\geq 0})^n\setminus\Delta(I)}$ is the simplex in $\R^n$ with vertices in the origin $A_0$
   and in the points $A_1,\dots,A_n$, where
   $$A_i=(\underbrace{0,\dots,0}_{\text{i-1}},\frac{s-(i-1)}{n-(i-1)},\underbrace{0,\dots,0}_{\text{n-i}}).$$
\end{theorem}

\section{Generic initial ideals}
   Let $S(n)=\K[x_1,\dots,x_n]$ be a polynomial ring over a field $\K$ of characteristic $0$.
   Let $\succ$ be the reversed lexicographical order on monomials in $S(n)$.
   Recall that $\succ$ is a total ordering defined as
   $$x_1^{p_1}\cdot\ldots\cdot x_n^{p_n}\;\succ\; x_1^{q_1}\cdot\ldots\cdot x_n^{q_n}$$
   if and only if, $\sum_{i=1}^np_i>\sum_{i=1}^n q_i$ or
   $\sum p_i=\sum q_i$ and there exists an index $k$ such that
   $p_\ell=q_\ell$ for all $\ell>k$ and $p_k< q_k$.

   For a homogeneous ideal $I\subset S(n)$, its \emph{initial ideal} $\ini(I)$
   is the ideal generated by leading terms
   of all elements of $I$. Recall that the \emph{leading term}
   $\ini(f)$ of a polynomial $f\in S$ is the greatest (with respect to
   the fixed order, here $\succ$) monomial summand of $f$. Initial ideals are of interest
   because they share many properties with the original ideals whereas
   they are easier to handle computationally, see e.g. \cite{HerHib11}.

   Even better behaved are \emph{generic initial ideals}. We describe now
   briefly how they are defined. To begin with, recall that $\GL(n,\K)$
   acts on $S(n)$ by the change of coordinates. The Borel subgroup $\T$
   of $\GL(n,\K)$ consists of upper triangular matrices
   $$\T=\left\{A\in\GL(n,\K):\; a_{ij}=0\;\mbox{ for all }\; j<i\right\}.$$
   Building upon ideas of Galligo \cite{Gal74}, Mark Green initiated in
   \cite{Gre98} a systematic study of generic initial ideals, see also \cite{GreSti98}.
   Theorem of
   Galligo \cite{Gal74} assures that for a homogeneous ideal $I$ and a generic
   choice of coordinates, the initial ideal $\ini(I)$ of $I$ is $\T$--fixed.
   In order to indicate that this property holds, we write $\gin(I)$ for $\ini(I)$.
   It follows from the same theorem that $\gin(I)$ is well defined, in particular
   uniquely determined by $I$.
   Generic initial ideals carry even more information on original ideal $I$
   than arbitrary initial ideals do. For example Hilbert functions of $I$
   and $\gin(I)$ are equal \cite[page 90]{GreSti98}.

\section{Limiting shape of a graded family of ideals}
   Following Mayes, we recall here a construction of a solid associated
   to a polynomial ideal with zero--dimensional support. To some
   extend this construction resembles that of Okounkov bodies \cite{LazMus09},
   even though the resulting body in our case need not to be convex.

   Let $\ibul=\left\{I_m\right\}$ be a graded family of ideals in $S(n)$, i.e.
   a collection of ideals $I_m\subset S(n)$ such that
   $$I_k\cdot I_\ell\subset I_{k+\ell}$$
   holds for all $k,\ell\geq 1$, see \cite[Definition 2.4.14]{PAG}.

   Assume that all ideals $I_m$ are of finite colength i.e.
   $\length(S(n)/I_m)$ is finite for all $m\geq 1$. To such a family
   one associates an invariant analogous to the volume of a graded
   linear series (this invariant is actually called a volume in \cite{May14},
   we prefer here however to follow the notation from \cite[Definition 2.3.40]{PAG}).
\begin{definition}[Multiplicity]
   The \emph{multiplicity} of a graded family $\ibul$ is the real number
   $$\mult(\ibul):=\limsup_{m\to\infty}\frac{\length(S(n)/I_m)}{m^n/n!}.$$
\end{definition}
   As customary we associate to a monomial $M=x_1^{m_1}\cdot\ldots\cdot x_n^{m_n}\in S(n)$
   a point $P(M)=(m_1,\ldots,m_n)$ in $\R^n$.
\begin{definition}[Newton polytope]
   Let $J$ be an arbitrary monomial ideal in $S(n)$. The \emph{Newton polytope}
   $P(J)$ of $J$ is the convex hull of the set
   $\left\{P(M)\in\R^n:\;M\in J\right\}$.
\end{definition}
   Given a graded family $\ibul$ of monomial ideals
   we are interested in the limiting shape of the Newton polytopes $P(I_m)$ defined as
   $$\Delta(\ibul):=\bigcup_{m=1}^{\infty}\frac1m P(I_m).$$
   For a monomial ideal $J\subset S(n)$, we denote by $Q(J)$ the closure of the complement of
   the Newton polytope $P(J)$ in the positive octant.
   It is well known, see e.g. \cite[Proposition 2.14]{May14}
   and \cite[Theorem 1.7 and Lemma 2.13]{Mus02} that
   if $J$ is supported on a zero-dimensional subscheme then
   \begin{equation}\label{eq:mult and vol}
      \mult(J)=n!\cdot \vol_{\R^n}(Q(J)).
   \end{equation}
   Here $\vol_{\R^n}$ denotes the standard Euclidean volume on $\R^n$
   normalized so that the unit cube $[0,1]^n$ has volume $1$.

   The asymptotic counterpart of sets $Q(I_m)$ for a graded
   family $\ibul$ of ideals is defined as
   $$\Gamma(\ibul)=\bigcap_{m=1}^{\infty} \frac1m Q(I_m).$$
   Of course $\Gamma(\ibul)$ is the closure of the complement of the convex set $\Delta(\ibul)$.
   Thus it is a coconvex body as studied in \cite{KhoTim13}.

   To an arbitrary ideal $J$, we associate a graded
   family of monomial ideals
   $$G_m:=\gin(J^{(m)}).$$
   We write then simply $\Gamma(J)$ for $\Gamma(\gbul)$.

\section{Star configurations}
   Let $s\geq n$ be a fixed integer and let
   $H_1,\ldots,H_s$ be general hyperplanes in $\P^n$.
   Each $n$ of these hyperplanes determine then a unique
   point in $\P^n$, the intersection point of all of them.
   Thus the set $Z$ consisting of all these points has
   exactly $\binom{s}{n}$ elements. Let $I=I_Z\subset S(n+1)$
   be the homogeneous ideal associated to the reduced
   scheme $Z$. We consider the graded family
   of symbolic powers of $I$
   $$I_m:=I^{(m)}.$$
   Symbolic powers of ideals supported on zero-dimensional
   reduced schemes are saturated, see e.g. \cite[Lemma 1.3]{ArsVat03}.
   A property of saturated ideals important for the construction
   carried out in the next section is the following result due to Green
   \cite[Theorem 2.21]{Gre98}.
\begin{proposition}[Green]
   Let $J$ be a saturated ideal in $S(n+1)$. Then
   no minimal generator of the generic initial
   ideal $\gin(J)$ contains the variable $x_{n+1}$.
\end{proposition}
   It follows that in the above set-up of a star configuration,
   even though the ideals $I^{(m)}$ are contained in $S(n+1)$,
   their generic initial ideals $\gin(I^{(m)})$ can be naturally considered
   as contained in $S(n)$. Thus their Newton polytopes
   and so also $\Delta(I)$ and $\Gamma(I)$
   are naturally contained in $\R^n$.

   Now we are in the position to proof Theorem \ref{thm:main points}.
\paragraph{Proof of Theorem \ref{thm:main points}.}
   Let $H_1,\dots,H_s$ be general hyperplanes in $\P^n$
   defined as zero-sets of linear polynomials $h_1\dots,h_s$.
   Let $Z$ be the set consisting of points where $n$
   of these hyperplanes meet. By generality $Z$ contains
   exactly $\binom{s}{n}$ mutually distinct points. Let $I$
   be the ideal of $Z$.

   It is enough to show that points $A_i$ are intersection points
   of the coordinate axes with the sets $\Delta(I)$ and $\Gamma(I)$.
   Indeed, taking this for granted for a while, we show how the
   assertion of the Theorem follows. To this end let $W$ denote
   the convex hull of the set $\left\{A_0, A_1,\dots,A_n\right\}$.
   Since $\Delta(I)$ is a convex set, its complement $\Gamma(I)$
   is for sure contained in $W$. On the other hand the volume
   of $W$ is
   $$\vol_{\R^n}(W)=\frac{1}{n!}\cdot\frac{s}{n}\cdot\frac{s-1}{n-1}\cdot\ldots\cdot\frac{s-(n-1)}{1}=\frac{1}{n!}\binom{s}{n}.$$
   According to \eqnref{eq:mult and vol} the volume of $\Gamma(I)$ is equal to
   $\frac{1}{n!}$ times the number of points supporting $I$. Hence $\vol_{\R^n}\Gamma(I)=\vol_{\R^n}(W)$
   and since both sets are closed and $\Gamma(I)\subset W$, they are equal.

   In fact, the above argument shows a little bit more.
   It is namely enough to prove that $A_1,\ldots,A_n$ are
   contained in $\Delta(I)$. This forces $\Gamma(I)$
   to be contained in $W$ and we are done again.

   Turning thus to the containment of points $A_1,\ldots,A_n$
   in $\Delta(I)$ by way of the warm up, we show that $A_1$
   appears in $P(\gin(I^{(n)}))$. To this end note first that
   the product $h_{1,2,\ldots,s}=h_1\cdot\ldots\cdot h_s$ is by construction
   an element of $I^{(n)}$. Since these forms are general, their product $h_{1,2,\ldots,s}$ contains the
   monomial $x_1^s$, which is then the leading term, so that
   $\gin(I^{(n)})$ contains this monomial as well. Thus, after scaling, we get the point
   $A_1=(\frac{s}{n},0,\dots,0)\in\frac1nP(\gin(I^{(n)}))$. We will give an additional interpretation
   of this vertex right after the current proof.

   \textit{Claim.} Now we will show that $\gin(I^{(n+1-k)})$ contains the monomial
   $x_{n+1-k}^{s+1-k}$.

   Since we work in characteristic zero and symbolic powers agree with differential powers
   via the Zariski--Nagata Theorem \cite[Section II.3.9]{Eis95}, the elements
   of degree $s+1-k$ in $I^{(n+1-k)}$ for $k=1,\ldots,n$ form a vector space. Its dimension is at least
   equal to $\binom{s}{k-1}$,
   the number of distinct products of $k-1$ forms $h_{j_1,\ldots,j_{k-1}}=h_{j_1}\cdot\ldots\cdot h_{j_{k-1}}$
   with $1\leq j_1<j_2<\ldots<j_{k-1}\leq s$.
   Indeed, all products of the form
   \begin{equation}\label{eq:form of elements}
      \frac{h_1\cdot\ldots\cdot h_s}{h_{j_1,\ldots,j_{k-1}}}
   \end{equation}
   vanish in all points of $Z$ to order at least $n+1-k$.

   \textit{Subclaim.}
   These products are moreover linearly independent. This can be seen as follows.
   Suppose that there is a linear combination with
   $$\sum_{1\leq j_1<\ldots<j_{k-1}\leq s}\lambda_{j_1,\ldots,j_{k-1}}\cdot \frac{h_1\cdot\ldots\cdot h_s}{h_{j_1,\ldots,j_{k-1}}}=0.$$
   Restricting this sum to the intersection $Y=H_{\ell_1}\cap\ldots\cap H_{\ell_{k-1}}$
   we obtain
   $$\lambda_{\ell_1,\ldots,\ell_{k-1}}\cdot h_{i_1,\ldots,i_{s+1-k}}=0$$
   for $\left\{i_1,\ldots,i_{s+1-k}\right\}=\left\{1,2,\ldots,s\right\}\setminus\left\{\ell_1,\ldots,\ell_{k-1}\right\}$.
   Since none of the forms $h_{i_1},\ldots,h_{i_{s+1-k}}$ vanishes identically
   on $Y$, it must be $\lambda_{\ell_1,\ldots,\ell_{k-1}}=0$. This argument
   works for arbitrary pick of $\ell_1,\ldots,\ell_{k-1}$ so that we can conclude
   that all $\lambda_{i_1,\ldots,i_{k-1}}$ vanish. This establishes the Subclaim.

   It is easy to see that there are $\binom{s}{k-1}-1$ monomials of degree $s+1-k$ which are
   strictly bigger than $x_{n+1-k}^{s+1-k}$ with respect to the order $\succ$. Hence
   there exists a non-zero linear combination $f$ of forms of degree $s+1-k$ as in \eqnref{eq:form of elements},
   with vanishing coefficients at all monomials strictly bigger than $x_{n+1-k}^{s+1-k}$ with respect
   to the order $\succ$. By the generality of the forms, the coefficient
   of $f$ at $x_{n+1-k}^{s+1-k}$ does not vanish. Hence $x_{n+1-k}^{s+1-k}$
   is the leading monomial of $f$ and thus it appears in $\gin(I^{n+1-k})$. This
   finishes the proof of Theorem \ref{thm:main points}.
\qed
\section{Some asymptotic invariants attached to a homogeneous ideal}
\subsection{$A_1$ and Waldschmidt constants}
   Now we give an interpretation of vertices $A_1$ and $A_n$.
   We recall first a couple of definitions.
\begin{definition}[Initial degree]
   Let $I$ be a homogeneous ideal in a polynomial ring. Then
   $$\mbox{the initial degree }\; \alpha(I)$$
   is defined as the minimal number $t$ such that there exists
   a non-zero element of degree $t$ in $I$.
\end{definition}
   The asymptotic counter-part of the initial ideal is called
   \emph{Waldschmidt constant} of $I$. It has been rediscovered
   recently and intensively studied by Harbourne, see e.g. \cite{BocHar10b}
   and \cite{GHVT13}.
\begin{definition}[Waldschmidt constant]
   The \emph{Waldschmidt constant} of $I$ is defined as
   $$\widehat{\alpha}(I):=\lim_{m\to\infty}\frac{\alpha(I^{(m)})}{m}.$$
\end{definition}
   The existence of the limit in the definition above follows from
   the sub-additivity of initial degrees. This also implies
   in a standard way that the Waldschmidt constant is in fact the
   infimum of the quotients $\frac{\alpha(I^{(m)})}{m}$.

   Waldschmidt
   constants for points forming a star configuration of $n$--wise
   intersection points of $s\geq n$ general hyperplanes in $\P^n$
   was computed in \cite[Example 8.3.4]{PSC}. We have
   $$\widehat{\alpha}(I)=\frac{s}{n}$$
   in this case. This corresponds to the non-zero coordinate of
   the vertex $A_1$.
\begin{remark}\label{rmk:Borel invariant}\rm
   In fact it is easy to see that $\alpha(I^{(m)})$ is the minimal coordinate
   of a point in $P(\gin(I^{(m)}))$ sitting on the $x_1$--axis. This follows from
   the Borel invariance of generic initial ideals. Indeed if $f$ is a monomial
   in $\gin(I^{(m)})$ divisible by some $x_i$, then also $x_1\cdot\frac{f}{x_i}$
   is contained in $\gin(I^{(m)})$. This implies that the minimal degree of a non-zero
   element in $I^{(m)}$ (which is of course equal to the degree of its leading term)
   can be read of the $x_1$--axis.
\end{remark}
\subsection{$A_n$ and asymptotic regularity}
   Now we want to interpret the last vertex $A_n$. We begin by recalling
   the notion of
   the Castelnuovo-Mumford regularity which is one of fundamental
   invariants determining the complexity of an ideal, see \cite[Section 20.5]{Eis95}
   and \cite[Section 1.8]{PAG}.
\begin{definition}[Castelnuovo-Mumford regularity]
   Let $J$ be a homogeneous ideal in a polynomial ring $S$ and let
   $$0\to\ldots\to\bigoplus_jS(-a_{ij})\to\ldots \bigoplus_jS(-a_{1j})\to\bigoplus_jS(-a_{0j})\to J\to 0$$
   be a minimal free resolution of $I$. Then the \emph{Castelnuovo-Mumford regularity} of $J$
   is the non-negative integer
   $$\reg(J):=\max_{i,j}\left\{\; a_{ij}-i\;\right\}.$$
\end{definition}
   In other words regularity of an ideal is governed by the maximal
   degree of generators of its syzygy modules. In particular
   we always have $\reg(J)\geq d(J)$, where $d(J)$ denotes the maximal degree
   in a minimal set of generators of $J$.
   For a Borel fixed ideal, Bayer and Stillman show in \cite[Proposition 2.9]{BaySti87}
   that one has actually the equality
   \begin{equation}\label{eq:reg for Borel inv}
      \reg(J)=d(J).
   \end{equation}
   Since the regularity is upper-semicontinuous in flat families
   one has always $\reg(J)\leq \reg(\ini(J))$
   and $\reg(J)=\reg(\gin(J))$ by \cite[Proposition 2.11]{BaySti87}.
   Now, the argument used in Remark \ref{rmk:Borel invariant} implies
   that $d(J)$ is detected by the generator of $J$ whose
   initial monomial is a power of $x_n$.

   Asymptotic regularity was studied by Cutkosky, Herzog and Trung in \cite{CHT99}.
   It is defined as the real number
   $$\areg(J)=\lim_{m\to\infty}\frac{\reg(J^m)}{m}.$$
   As remarked
   in \cite[Theorem 1.1]{CHT99} it is equal to $\lim_{m\to\infty}\frac{d(J^m)}{m}$.

   Motivated more by geometry Cutkosky, Ein and Lazarsfeld study in \cite{CEL01}
   a new invariant associated to an ideal sheaf. Their definition works for ideal
   sheaves on arbitrary projective varieties, we prefer however to restrict our
   attention here to the projective space setting.
\begin{definition}[the $s$--invariant]
   Let $\cali$ be an ideal sheaf in $\calo_{\P^n}$ and let
   $\pi:X\to\P^n$ be the blow up of $\cali$ with exceptional
   divisor $F$. The \emph{$s$--invariant} of $\cali$ is defined as
   $$s(\cali):=\min\left\{t\in\R:\; t\cdot\pi^*(\calo_{\P^n}(1))-F\;\mbox{ is nef}\right\}.$$
\end{definition}
   If $\cali$ is the sheafification of a homogeneous ideal $I$, then we
   have
   \begin{equation}\label{eq:areg = s}
      \areg(I)=s(\cali),
   \end{equation}
   see \cite[Theorem 5.4.22]{PAG}.

   Asymptotic regularity of symbolic powers was studied recently by
   Cutkosky and Kurano in \cite{CutKur11}. It is defined as
   $$\asreg(I):=\lim_{m\to\infty}\frac{\reg(I^{(m)})}{m}.$$
   It follows immediately
   from \cite[Theorem 4.6]{CutKur11} that also
   \begin{equation}\label{eq:asreg = s}
      \asreg(I)=s(\cali).
   \end{equation}

   Since generic initial ideal $\gin(I)$ of a homogeneous ideal $I$ in $S(n+1)$
   contains among generators a power of $x_n$, see \cite[Corollary 2.9]{May14},
   combining \eqnref{eq:reg for Borel inv} with \eqnref{eq:areg = s} and \eqnref{eq:asreg = s}
   we get that the non-zero coordinate of $A_n$ equals
   $$\areg(I)=\asreg(I)=s(\cali).$$
\begin{remark}
   It would be interesting to know if the remaining vertices $A_2,\ldots,A_{n-1}$
   can be interpreted along similar lines.
\end{remark}
   We conclude this note with an example showing that one cannot reverse
   Theorem \ref{thm:main points}, i.e. the shape of $\Gamma(I)$ does
   not determine $I$ to be the ideal of a star configuration of points.
\begin{example}[Points with star configuration like limiting shape]\rm
   Let
   $I$ be the ideal of intersection points $P_1,\ldots,P_6$ of a smooth conic $C$ with
   a general cubic curve $D$. Then obviously
   $\widehat{\alpha}(I)=2$ and $\areg(I)=3$. Hence
   $\Gamma(I)$ is the triangle with vertices $A_0=(0,0)$, $A_1=(2,0)$ and $A_2=(0,3)$
   but the points $P_1,\ldots,P_6$ do not form a star configuration.
\end{example}

%\paragraph*{\emph{Acknowledgement.}}

%*****************************************************************************

%***************************************************************************** % Addresses

\bigskip \small

\bigskip
   Marcin Dumnicki,
   Jagiellonian University, Institute of Mathematics, {\L}ojasiewicza 6, PL-30-348 Krak\'ow, Poland

\nopagebreak
   \textit{E-mail address:} \texttt{Marcin.Dumnicki@im.uj.edu.pl}

\bigskip
   Tomasz Szemberg,
   Instytut Matematyki UP,
   Podchor\c a\.zych 2,
   PL-30-084 Krak\'ow, Poland.

\nopagebreak
   \textit{E-mail address:} \texttt{szemberg@up.krakow.pl}

\bigskip
   Justyna Szpond,
   Instytut Matematyki UP,
   Podchor\c a\.zych 2,
   PL-30-084 Krak\'ow, Poland.

\nopagebreak
   \textit{E-mail address:} \texttt{szpond@up.krakow.pl}

\bigskip
   Halszka Tutaj-Gasi\'nska,
   Jagiellonian University, Institute of Mathematics, {\L}ojasiewicza 6, PL-30-348 Krak\'ow, Poland

\nopagebreak
   \textit{E-mail address:} \texttt{Halszka.Tutaj@im.uj.edu.pl}

%*****************************************************************************

\end{document}